\documentclass{article}
\usepackage{amsmath,amssymb,amsthm,bm}
\usepackage[dvipdfmx]{graphicx}
\usepackage{xcolor}
\usepackage{hyperref}

\title{Existence of Really Perverse Central Configurations in the Spatial $N$-Body Problem}
\author{Mitsuru Shibayama}
\date{}

\newtheorem{theorem}{Theorem}
\newtheorem{definition}{Definition}

\begin{document}
\maketitle

\begin{abstract}
We construct explicit examples of really perverse central configurations in the spatial Newtonian $N$-body problem.
A central configuration is called really perverse if it satisfies the central configuration equations
for two distinct mass distributions having the same total mass and the same center of mass.
While such configurations were previously known only in the planar case for large values of $N$,
we prove the existence of spatial really perverse central configurations for $N=27,\dots,55$.
\end{abstract}

\section{Introduction}

Central configurations play a fundamental role in celestial mechanics.
They generate self--similar solutions of the Newtonian $N$--body problem,
and they appear naturally in the study of total collisions, relative equilibria, and
bifurcations of energy manifolds.

Chenciner introduced the notion of \emph{perverse solutions}
and \emph{really perverse solutions} of the planar $n$--body problem.
A solution is called \emph{perverse} if it solves the equations for two distinct mass distributions,
and \emph{really perverse} if moreover the total mass and the center of mass remain unchanged.

This notion is strongly motivated by the theory of choreographies.
As pointed out by Chenciner \cite{Chenciner04}, if a planar choreography exists with non--equal masses,
then by replacing each mass with the mean mass $M/n$ one obtains another admissible mass distribution,
which automatically makes the choreography really perverse.

The planar case has been investigated.
For $n=2$ no perverse solution exists, since the mass distribution is uniquely determined.
For $n=3$, Albouy and Moeckel \cite{AM} proved that the only perverse solutions are collinear homographic motions,
and that no really perverse solution exists in this case.
For $n=4$, Chenciner \cite{Chenciner03} showed that any perverse solution must satisfy strong geometric constraints,
and concluded that no really perverse planar solution exists for $n\le4$. 

The first examples of really perverse solutions in the planar problem were later found
by Chenciner--Sim\'o for configurations consisting of a central mass surrounded by
 homothetic regular polygons, with the smallest example occurring for $n=1369$
bodies (see the last part of \cite{Chenciner03}).

More recently, Bang, Chenciner and Sim\'o \cite{BangChencinerSimo2025} refined this construction and obtained families
of planar really perverse solutions with significantly smaller numbers of bodies,
down to $N=474$ in the most efficient case.

In contrast to the planar situation, the spatial case has remained completely open.
No example of a really perverse solution or central configuration was known in dimension three.
The purpose of the present paper is to fill this gap by constructing explicit examples of
\emph{spatial} really perverse central configurations.

\section{Main Result and its Proof}

Consider the Newtonian $N$-body problem in $\mathbb{R}^3$: 
\begin{equation}\label{eq:newton}
\ddot Q_k = -\sum_{j\neq k}\frac{M_j(Q_k-Q_j)}{|Q_k-Q_j|^3},
\qquad k=1,\dots,N.
\end{equation}

A configuration $\bm{c}=(c_1,\dots,c_N)\in(\mathbb{R}^3)^N$ is a central configuration if there exists $\lambda\in\mathbb{R}$ such that
\begin{equation}\label{eq:cc}
\lambda M_k c_k = -\sum_{j\neq k}\frac{M_kM_j(c_k-c_j)}{|c_k-c_j|^3}.
\end{equation}

We introduce the following definition:

\begin{definition}
A central configuration $\bm{c}$ is called \emph{really perverse}
if it satisfies \eqref{eq:cc} for two distinct mass distributions
$(M_1,\dots,M_N)\neq(M_1',\dots,M_N')$
with the same total mass and the same center of mass.
\end{definition}

If $\bm{c}$ is a really perverse central configuration, then the associated homothetic motion is a really perverse solution.

\begin{theorem}\label{thm:main}
There exist really perverse central configurations of the spatial Newtonian $N$-body problem
for all $N=27,28,\dots,55$.
\end{theorem}

Fix $n \ge 2$ and consider the symmetric configuration
\begin{align}\label{eqn:symsol}
q_0 = (0, 0, 0), ~
q_k = r(t) (\cos \tfrac{2\pi  k}{n},\sin \tfrac{2\pi  k}{n},  0)  ~ (k=1, \dots, n), 
~ q_{\pm} = r(t)(0, 0, \pm \alpha)
\end{align}
with masses $m_0,m_1,m_2$ respectively.

\begin{figure}[h]
\begin{center}
\includegraphics[width=2in]{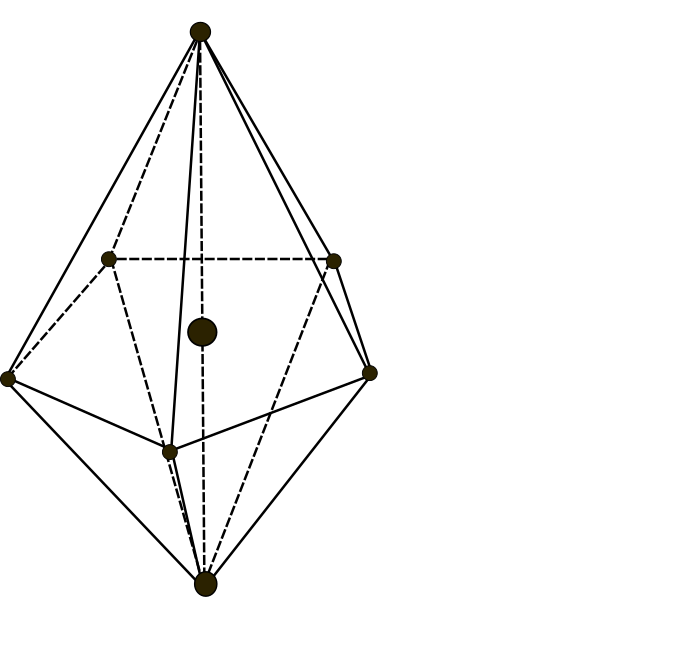}
\caption{The symmetric configuration \eqref{eqn:symsol}}
\end{center}
\end{figure}

The equations of motion take the form
\begin{align*}
\ddot{q}_k &= -\frac{m_0}{|q_k|^3}q_k 
- \sum_{j \neq k} \frac{m_1(q_k - q_j)}{|q_k - q_j|^3}
- \sum \frac{m_2(q_k - q_\pm)}{|q_k - q_\pm|^3}, \\
\ddot{q}_\pm &= -\frac{m_0}{|q_\pm|^3}q_\pm 
- \sum_j \frac{m_j(q_\pm - q_j)}{|q_\pm - q_j|^3}
- \frac{m_2(q_+ - q_-)}{|q_+ - q_-|^3}.
\end{align*}
We note that the center of mass is at the origin.

Substituting \eqref{eqn:symsol} yields
\begin{align*}
\ddot{r} &= -\frac{m_0}{r^2}
- \sum_{j \neq 0} \frac{m_1(r - e^{2 \pi i j/n}r)}{|r - e^{2 \pi i j/n}r|^3}
- \sum \frac{m_2 r(1, \pm  \alpha)}{|r(1, \pm  \alpha)|^3}, \\
 \alpha\ddot{r} &= -\frac{m_0}{ \alpha^2 r^2}
- \sum_j \frac{m_j(-r e^{2\pi i j/n},  \alpha r)}{|(-r e^{2\pi i j/n},  \alpha r)|^3}
- \frac{m_2(2  \alpha r)}{|2  \alpha r|^3}.
\end{align*}

Assume $\ddot r=-1/r^2$. Then the equations reduce to
\begin{align}
1 &= m_0 + m_1 H_n + \frac{2m_2}{(1+\alpha^2)^{3/2}}, \label{eq:A}\\
\alpha^3 &= m_0 + \frac{n\alpha^3m_1}{(1+\alpha^2)^{3/2}} + \frac{m_2}{4}, \label{eq:B}\\
M &= m_0 + nm_1 + 2m_2 \label{eq:C}
\end{align}
where
\[H_n = \sum_{j=1}^{n-1} \frac{1 - e^{2\pi ij/n}}{|1 - e^{2\pi ij/n}|^3}
= \sum_{j=1}^{n-1} \frac{1 - \cos(2\pi j/n)}{|1 - e^{2\pi ij/n}|^3}
= \sum_{j=1}^{n-1} \frac{1}{4|\sin(\pi j/n)|}.\]

We also impose
\[
M = m_0 + n m_1 + 2m_2.
\]

Eliminating $m_0$ from \eqref{eq:A}--\eqref{eq:C}, we obtain
\begin{align}
M-1 &= nm_1\Bigl(1-\frac{H_n}{n}\Bigr)
+2\Bigl(1-\frac{1}{(1+\alpha^2)^{3/2}}\Bigr)m_2,\\
M-\alpha^3 &= nm_1\Bigl(1-\frac{\alpha^3}{(1+\alpha^2)^{3/2}}\Bigr)
+\frac{7}{4}m_2.
\end{align}

The determinant of the coefficient matrix is
\[
f_n(\alpha)=
\frac74\Bigl(1-\frac{H_n}{n}\Bigr)
-2\Bigl(1-\frac{2}{(1+\alpha^2)^{3/2}}\Bigr)
\Bigl(1-\frac{\alpha^3}{(1+\alpha^2)^{3/2}}\Bigr).
\]

\begin{figure}[h]
\begin{center}
\includegraphics[width=2in]{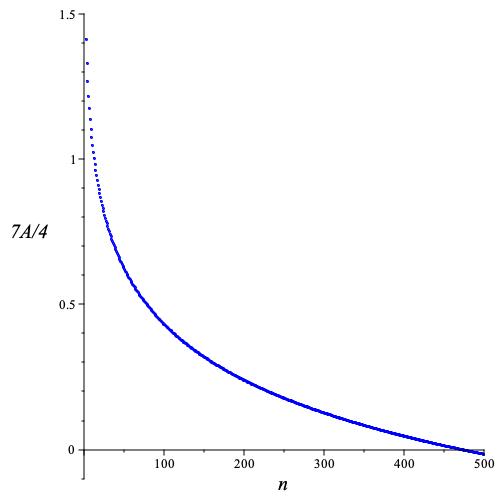}
\includegraphics[width=2in]{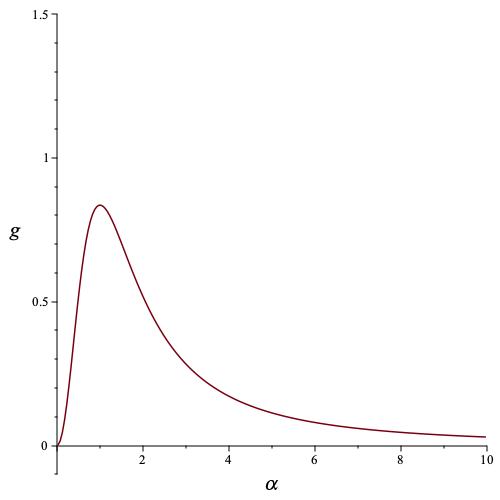}
\caption{Graphs of $\frac74\Bigl(1-\frac{H_n}{n}\Bigr)$ and $2\Bigl(1-\frac{2}{(1+\alpha^2)^{3/2}}\Bigr)
\Bigl(1-\frac{\alpha^3}{(1+\alpha^2)^{3/2}}\Bigr)$.}
\end{center}
\end{figure}

We first choose $\alpha>0$ so that $f_n(\alpha)=0$.
Let
\[
g(\alpha)
= 2\Bigl(1-\frac{2}{(1+\alpha^2)^{3/2}}\Bigr)
\Bigl(1-\frac{\alpha^3}{(1+\alpha^2)^{3/2}}\Bigr).
\]
Then $g$ attains its maximum value $\frac94-\sqrt{2}$ at $\alpha=1$ and tends to $+0$ as $\alpha \to +\infty$.

On the other hand, the quantity
\[
\frac74\Bigl(1-\frac{H_n}{n}\Bigr)\]
is strictly decreasing in $n$ and satisfies
\[
\frac74\Bigl(1-\frac{H_n}{n}\Bigr) < \frac94-\sqrt{2}
\qquad \text{for all } n\ge24.
\]
Let 
\[ A = 1- \frac{H_n}{n}. \] 
Moeckel and Sim\'o  \cite{MoeckelSimo} showed 
 \[ A >0 \quad (n \le 472), \qquad A<0 \quad (n \ge 473). \]
Hence, by the intermediate value theorem, 
for $24 \le n \le 472$, 
there exists $\alpha_0 \in (0, 1)$
such that $f_n(\alpha_0)=0$.

Next we determine $M$ so that the first two equations of the system
are equivalent. Let 
\[  B:= 1- \frac{\alpha_0^3}{(1+\alpha_0^2)^{3/2}}, \quad C=1-\frac{1}{(1+\alpha_0^2)^{3/2}} >0.\]
We determine $M$ such that 
 \[ \frac{M-1}{M-\alpha_0^3 } =\frac{A}{B}, 
 \]
 that is 
 \begin{align*}
M&=
\frac{B-\alpha_0^3 A}{B-A}.
 \end{align*}
 Let $m_2=u$ be a parameter, and 
 determine 
$m_1$ and $m_0$ by 
 \begin{align*}
 m_1 &= \frac{1}{n A} (M-1 -2C m_2)= \frac{1}{n A} (M-1 -2C u), \\
 m_0 &=M - n m_1 -2 m_2 
  =(C-A)(  \frac{\alpha^3}{ B-A}  +\frac{2}{A}u).
 \end{align*}

 
Since 
 \[ f(\alpha_0)= \frac{7}{4}A -2 BC=0, \quad C <1-\frac{1}{2^{3/2}},  \]
 we obtain
 \[ \frac{B}{A}= \frac{7}{8C}> \frac{7}{8}  (1-\frac{1}{2^{3/2}})^{-1} \approx 1.353553390593274>1. \]
 Therefore
 $B>A$ and 
 hence $M>1$.
 For 
 \[ 0<u < \frac{M-1}{2C}=:u_0,  \]
  $m_1$ is positive.
 
If $C<A$, $m_0$ is also positive. 
The inequality is difficult to prove analytically.
We verify this inequality using rigorous interval arithmetic (INTLAB).
 See Table \ref{tb:cap}.
 For example, 
 $C-A =0.128657457067_{401}^{413}$ represents 
 \[ C-A \in [0.128657457067401, 0.128657457067413].\]

For $n = 24, \dots, 52$, we rigorously verify that $C - A > 0$, which yields positive solutions
$(m_0,m_1,m_2)$  of the system for  $u \in (0, u_0)$.
Therefore, really perverse central configurations exist for these values of $n$.
This completes the proof of Theorem~\ref{thm:main}.

\begin{table}[htbp]
  \centering
\begin{tabular}{|c|c|c|c|c|c|}
\hline
$n$ & $\alpha_0$ & $M$ & $C-A$ & $u_0$  \\\hline
$24$ & $0.922477314112_{091}^{111}$ & $1.47638514245_{6963}^{7427}$ & $0.128657457067_{401}^{413}$ & $0.3950798816_{59850}^{60244}$ \\ 
$25$ & $0.8702462710626_{66}^{87}$ & $1.63954060886_{3972}^{4361}$ & $0.103000083155_{111}^{124}$ & $0.560273281629_{642}^{996}$ \\
$26$ & $0.8365637482706_{20}^{40}$ & $1.697355545841_{204}^{543}$ & $0.087284996115_{715}^{729}$ & $0.635368964774_{340}^{665}$ \\
$27$ & $0.8105444160356_{63}^{83}$ & $1.722173365640_{276}^{582}$ & $0.075675469693_{770}^{784}$ & $0.679805754954_{832}^{138}$ \\
$28$ & $0.7890092982849_{51}^{71}$ & $1.731857491632_{300}^{577}$ & $0.066456471337_{093}^{108}$ & $0.708954889882_{619}^{908}$ \\
$29$ & $0.770499013641_{090}^{110}$ & $1.73346557503_{9926}^{0182}$ & $0.058837688928_{116}^{131}$ & $0.729168230930_{496}^{773}$ \\
$30$ & $0.7541988921012_{54}^{74}$ & $1.730407286581_{044}^{281}$ & $0.052377057360_{983}^{999}$ & $0.743661093181_{250}^{515}$ \\
$31$ & $0.7395993040349_{69}^{89}$ & $1.724532530152_{108}^{330}$ & $0.046797902783_{883}^{899}$ & $0.754266535491_{748}^{003}$ \\
$32$ & $0.7263564825238_{58}^{78}$ & $1.716923788961_{257}^{464}$ & $0.041913757083_{905}^{921}$ & $0.762115803421_{607}^{853}$ \\
$33$ & $0.7142257236271_{72}^{92}$ & $1.708248672037_{514}^{709}$ & $0.037592155281_{122}^{138}$ & $0.767947999932_{456}^{695}$ \\
$34$ & $0.70302565585_{6998}^{7018}$ & $1.698934776483_{611}^{796}$ & $0.033735214280_{357}^{373}$ & $0.772267320016_{562}^{794}$ \\
$35$ & $0.6926175469625_{71}^{91}$ & $1.689263581439_{451}^{627}$ & $0.030268355827_{115}^{131}$ & $0.775429564837_{473}^{699}$ \\
$36$ & $0.6828925650940_{75}^{95}$ & $1.679424017592_{502}^{669}$ & $0.027133348724_{076}^{093}$ & $0.777692776313_{726}^{947}$ \\
$37$ & $0.6737635473245_{11}^{31}$ & $1.669544528846_{037}^{196}$ & $0.024283804520_{173}^{189}$ & $0.779248363506_{408}^{624}$ \\
$38$ & $0.6651594666655_{37}^{57}$ & $1.65971301848_{5976}^{6128}$ & $0.021682145366_{774}^{791}$ & $0.780241018463_{549}^{760}$ \\
$39$ & $0.6570215890663_{66}^{86}$ & $1.649989659683_{357}^{503}$ & $0.019297495732_{172}^{189}$ & $0.780781890895_{911}^{118}$ \\
$40$ & $0.6493007293582_{18}^{38}$ & $1.640415344753_{219}^{360}$ & $0.017104175967_{791}^{808}$ & $0.780957552642_{056}^{260}$ \\
$41$ & $0.6419552449709_{06}^{26}$ & $1.631017384269_{305}^{440}$ & $0.015080600564_{736}^{753}$ & $0.780836246655_{377}^{576}$ \\
$42$ & $0.634949538711_{199}^{219}$ & $1.621813425081_{826}^{955}$ & $0.013208456019_{241}^{258}$ & $0.780472335423_{192}^{388}$ \\
$43$ & $0.6282529212561_{78}^{98}$ & $1.612814187984_{423}^{548}$ & $0.011472076490_{259}^{276}$ & $0.779909526377_{221}^{413}$ \\
$44$ & $0.6218387331925_{32}^{52}$ & $1.604025407352_{340}^{460}$ & $0.009857962288_{385}^{402}$ & $0.779183248840_{005}^{195}$ \\
$45$ & $0.6156836578161_{28}^{48}$ & $1.595449221719_{335}^{452}$ & $0.008354403402_{385}^{402}$ & $0.778322431220_{738}^{926}$ \\
$46$ & $0.6097671764602_{31}^{52}$ & $1.587085180738_{135}^{248}$ & $0.006951181530_{579}^{596}$ & $0.777350847125_{508}^{693}$ \\
$47$ & $0.6040711318967_{04}^{24}$ & $1.578930980460_{547}^{656}$ & $0.005639331642_{870}^{887}$ & $0.776288146931_{273}^{455}$ \\
$48$ & $0.5985793747799_{68}^{88}$ & $1.570983003894_{489}^{594}$ & $0.004410949277_{764}^{781}$ & $0.775150656732_{715}^{896}$ \\
$49$ & $0.593277474673_{398}^{418}$ & $1.563236720512_{177}^{281}$ & $0.003259033392_{841}^{858}$ & $0.773952003112_{748}^{927}$ \\
$50$ & $0.5881524818550_{31}^{51}$ & $1.555686982631_{204}^{305}$ & $0.002177357151_{937}^{955}$ & $0.772703606030_{756}^{934}$ \\
$51$ & $0.5831927294515_{33}^{53}$ & $1.548328245774_{189}^{286}$ & $0.001160360880_{167}^{184}$ & $0.771415070822_{446}^{622}$ \\
$52$ & $0.5783876678962_{40}^{60}$ & $1.541154732586_{232}^{327}$ & $0.000203062767_{916}^{934}$ & $0.770094502289_{658}^{832}$ \\
$53$ & $0.57372772551_{5997}^{6017}$ & $1.534160554587_{635}^{726}$ & $-0.000699016096_{274}^{257}$ & $0.768748758099_{243}^{414}$ \\
\hline
\end{tabular}
  \caption{Rigorous bounds obtained by interval arithmetic (INTLAB)}
  \label{tb:cap}
\end{table}


\begin{thebibliography}{9}
\bibitem{AM}
A. Albouy and R.~B. Moeckel, The inverse problem for collinear central configurations, Celestial Mech. Dynam. Astronom. {\bf 77} (2000), 77--91 (2001)

\bibitem{BangChencinerSimo2025}
D. Bang, A. Chenciner and C. Sim\'o~i~Torres, Really perverse periodic solutions of the planar $N$-body problem, Celestial Mech. Dynam. Astronom. {\bf 137} (2025), Paper No. 26, 15 pp.

\bibitem{Chenciner03}
A. Chenciner, Perverse solutions of the planar $n$-body problem, Ast\'erisque No. 286 (2003), 
249--256

\bibitem{Chenciner04}
A. Chenciner, Are there perverse choreographies?, in {\it New advances in celestial mechanics and Hamiltonian systems}, 63--76, Kluwer/Plenum, New York, 

\bibitem{MoeckelSimo}
R.  Moeckel and C. Sim\'o, Bifurcation of spatial central configurations from planar ones, {\it SIAM J. Math. Anal.} {\bf 26} (1995),  978--998
\end{thebibliography}
\end{document}